\documentclass[a4paper,12pt]{article}
\usepackage{amssymb,amsmath,amsfonts}
\usepackage{bm}
\usepackage{cite}

\def\beq{\begin{equation}}
\def\eeq{\end{equation}}
\def\bea{\begin{eqnarray}}
\def\eea{\end{eqnarray}}
\def\nn{\nonumber}
\def\bra#1{\left\langle #1\right|}
\def\ket#1{\left| #1\right\rangle}
\def\bracket#1#2{\left\langle #1 \left| \right.\! #2 \right\rangle}

\def\U{{\cal U}}
\def\A{{\cal A}}
\def\hf{\frac{1}{2}}

\def\qn#1{\left[#1\right]_q}

\def\bn1#1#2{ \left[ \begin{array}{c} #1 \\ #2 \end{array} \right]_q }
\def\BH#1#2{ {}_{#1}\phi_{#2}}
\def\BHA#1#2#3#4{ \left[ \begin{array}{c} #1 \\ #2 \end{array}; \,#3 \,; #4 \right] }
\def\Haar#1{ H\big[ \; #1 \; \big]}
\renewcommand{\theequation}{\arabic{section}.\arabic{equation}}

\begin{document}
\thispagestyle{empty}

\vspace*{3cm}

\begin{center}
{\LARGE\sf
  Coherent State on \bm{$SU_q(2)$} Homogeneous Space
}

\bigskip\bigskip
N. Aizawa

\bigskip
\textit{
Department of Mathematics and Information Sciences, \\
Graduate School of Science, 
Osaka Prefecture University, \\
Nakamozu Campus, Sakai, Osaka 599-8531, Japan}\\

\bigskip
and

\bigskip
R. Chakrabarti \\

\bigskip
\textit{
Department of Theoretical Physics, \\
University of Madras, \\
Guindy Campus, Chennai 600 025, India
}

\end{center}

\vfill
\begin{abstract}
 The generalized coherent states for quantum groups introduced by 
Jur\v{c}o and \v{S}\v{t}ov\'{i}\v{c}ek are studied 
for the simplest example $ SU_q(2) $ in full detail. It is shown that the 
normalized $ SU_q(2) $ coherent states enjoy the property of 
completeness, and allow a resolution of the unity. This feature is expected to
play a key role in application of these coherent states in physical models. 
The homogeneous space of $ SU_q(2),$ $i.e.$ the  $q$-sphere of Podle\'s, is 
reproduced in complex coordinates by using the coherent states. Differential 
calculus in the complex form on the homogeneous space is developed. 
High spin limit  of the $ SU_q(2) $ coherent states is also 
discussed. 
\end{abstract}

%%%%%%%%%%%%%%%%%%%%%%%%%%%%%%%%%%%%%%%%%%%%%%%%%%%%%%%%%%%%%
%
%     Introduction
%
%%%%%%%%%%%%%%%%%%%%%%%%%%%%%%%%%%%%%%%%%%%%%%%%%%%%%%%%%%%%
%
\clearpage
\setcounter{page}{1}
\section{Introduction}
\label{intor}

  Generalization of boson coherent state based on 
a group theoretical viewpoint allows us to define 
coherent states for arbitrary Lie groups \cite{Per1,ACGT}. 
Extensive study of such generalized coherent states has revealed 
physical and mathematical richness of the notion. 
Physical applications \cite{Per, KS, ZFG} of the generalized coherent states 
are found in many instances such as quantum optics, semiclassical quantization 
of systems with spin degrees of freedom, transition from a pure state to mixed 
state dynamics during a nuclear collision, and so on.  
Mathematically, they provide \cite{Per, KS, ZFG} a natural framework to 
study geometric structure of the homogeneous space of 
the group under consideration. 

  On the other hand, advent of quantum groups and quantum 
algebras \cite{D1,D2,J1,J2,FRT,Ma,Wor} popularized deformations of Lie groups 
and Lie algebras in theoretical physics and mathematics. 
It is known that many properties 
of Lie groups and Lie algebras have their deformed counterparts.  
Especially, similarity in representation theories is remarkable. 
This motivates many researchers to consider coherent state for 
quantum algebras. Here we mention the pioneering works \cite{ArCo, Bie},
and attempts \cite{Que,Yu,Jur,MTZ,APVV} to construct the generalized coherent 
states of the type \cite{Per1,ACGT} in the context of quantum groups. 
The generalized coherent state is a vector in the representation space of 
the Lie group, and it is regarded as a function on homogeneous space 
with respect to a fiducial vector. For many Lie groups, such coherent states 
are obtained employing the duality between the group and its algebra. 
To see this, it may be enough to recall that the representation of a 
Lie group used in the construction of coherent states is obtained by the 
exponential mapping from a Lie algebra to the corresponding Lie group. 
Therefore generalized coherent states for quantum groups 
should embody the duality of a quantum group and its
quantum algebra. Such a definition was introduced by 
Jur\v{c}o and \v{S}\v{t}ov\'{i}\v{c}ek \cite{JS}. They developed 
a general theory of the quantum group version of generalized coherent states 
starting with the $q$-analogue of the Iwasawa decomposition and 
using the quantum double technique. 
The coherent states in  \cite{JS} reflect full Hopf algebra structure 
of quantum groups so that they may be the most plausible generalization of 
the coherent states of \cite{Per1,ACGT}. 
Nevertheless, to our knowledge, few works have been done on the coherent 
states for the quantum groups till today.

  In the present work, taking the simplest quantum group $ SU_q(2) $ as 
an example, we study the coherent states defined in  \cite{JS} in full detail.
Our construction of the coherent states is basically same as in \cite{JS}, 
but there is a slight difference. A central object for the coherent state 
construction is the universal $\cal T$-matrix, referred to as the canonical 
element in  \cite{JS}. Instead of using the Iwasawa decomposition as done in
\cite{JS} we derive the universal $\cal T$-matrix via the method of 
Fronsdal and Galindo \cite{FG}, since it makes clear that the universal 
$\cal T$-matrix is nothing but a quantum analogue of the exponential mapping 
from a Lie algebra to a Lie group. Anticipating that the $ SU_q(2) $ coherent 
states may have applications in developing field theories on noncommutative 
spaces, we put emphasis on the resolution of the unit operator and the complex 
description of the $ SU_q(2) $ homogeneous space, $i.e.,$ the Podle\'s 
\cite{Pod} $q$-sphere. Noncommutativity of the variables parametrizing the 
coherent states for quantum groups was discussed earlier, for instance, 
in \cite{MS}.

  The plan of this paper is as follows: The next preliminary section starts 
reviewing the quantum algebra $U_q[su(2)]$, and then the quantum group 
$ SU_q(2) $ is introduced as a dual algebra in the sense of \cite{FG}. The 
$*$-structure of $ SU_q(2) $ is studied in some detail emphasizing the 
$*$-conjugation properties of the noncommutative generators of the function 
algebra. We construct the $ SU_q(2) $ coherent states  and study their 
properties in \S \ref{SUq2CS}. Especially, the resolution of unity is proved in 
an algebraic setting. In \S \ref{CSqsphere} we show that the $ SU_q(2) $ 
coherent states naturally give the complex description of $q$-sphere of 
Podle\'s. Differential calculus on the $q$-sphere in complex description is 
explicitly provided. The Fock-Bargmann type representations of 
$ U_q[su(2)] $ is considered in \S \ref{CSrepU} as an application of the 
resolution of unity. High spin limit of the $ SU_q(2) $ coherent state 
is investigated in \S \ref{Contra-and-Ent}. It is shown that the limit gives a 
contraction to the coherent state of a quantum Heisenberg group. 
%In \S \ref{Entangled} we show that a construction of the entangled 
%two particle coherent states for the classical $ SU(2)$ group also applies 
%to the quantum group setting. Employing this we obtain the entangled 
%$SU_q(2)$ coherent states. 
We devote \S \ref{CRem} to concluding remarks.

%%%%%%%%%%%%%%%%%%%%%%%%%%%%%%%%%%%%%%%%%%%%%%%%%%%%%%%%%%%%%
%
%     Duality and $*$-structure
%
%%%%%%%%%%%%%%%%%%%%%%%%%%%%%%%%%%%%%%%%%%%%%%%%%%%%%%%%%%%%
%

\section{Duality and $*$-structure}
\label{Pre}

 The well-known \cite{D1,D2,J1,J2} quantum algebra $ \U = U_q[su(2)] $
endowed with a Hopf $*$-structure is 
generated by three elements $ J_{\pm}, \ J_0 $ subject to the relations 
\beq
  [J_0, J_{\pm} ] = \pm J_{\pm}, \qquad
  [J_+, J_- ] = [2J_0]_q,
  \label{Uqdef}
\eeq
where the $q$-deformed construct $[{\cal X}]$ reads:
\beq
   [{\cal X}]_q = \frac{ q^{\cal X} - q^{- {\cal X}} }
{ q - q^{-1} }. \label{qnum1}
\eeq
The coproduct $ \Delta$, 
the counit $ \epsilon $ and the antipode $S$ maps  given below
\bea
  & & \Delta(J_0) = J_0 \otimes 1 + 1 \otimes J_0, \qquad
      \Delta(J_{\pm}) = J_{\pm} \otimes q^{J_0} + q^{-J_0} \otimes J_{\pm},
      \label{Uqcopro} \\[3pt]
  & & \epsilon( J_0) = \epsilon( J_{\pm} ) = 0, 
      \label{Uqcounit} \\[3pt]
  & & S(J_0) = - J_0, \qquad S(J_{\pm}) = - q^{\pm 1} J_{\pm},
      \label{Uqanti}
\eea
together with the $*$-involution for $ q \in {\mathbb R} $ 
\beq
  J_+^* = J_-, \qquad J_-^* = J_+, \qquad J_0^* = J_0
  \label{star}
\eeq
satisfy the axioms of the Hopf $*$-algebra. A set of monomials 
\beq
E_{k\ell m} = J_+^k\, J_0^{\ell}\, J_-^m,
\qquad k, \ell, m \in {\mathbb Z}_{\geq 0}
\label{Ubasis}
\eeq
provide the basis of the universal enveloping algebra $ \U. $

  The algebra dual to $ \U $ is the quantum group $ \A = SU_q(2). $ One can 
determine the basis elements of the algebra $\A$ and their Hopf structure 
\textit{\`a la} Fronsdal and Galindo \cite{FG} where they started from two 
parameter deformation of $ SL(2), $ and reached its dual algebra $\U.$ 
Reversing their construction, we derive the basis
vectors $ e^{k \ell m} $ of the quantum group $\A$ by employing their dual 
relations with the known basis set $ E_{k \ell m} $ of $\U$:  
\beq
   \langle \; e^{k\ell m}, \; E_{k'\ell'm'} \, \rangle 
   = \delta^k_{k'}\, \delta^{\ell}_{\ell'}\, \delta^m_{m'}.
   \label{paring1}
\eeq
The two sets of structure constants defined below express the duality of 
$ \U $ and $ \A $ as follows:
\bea
  & & E_{k\ell m} E_{k' \ell' m'} = \sum_{abc} f_{k\ell m \ k'\ell'm'}^{\quad abc} 
      E_{abc},
     \quad 
      \Delta(E_{k\ell m}) = \sum_{abc \atop a'b'c'} g^{abc \ a'b'c'}_{\ \  k\ell m} 
      E_{abc} \otimes E_{a'b'c'},
      \label{UqStc} \\
  & & e^{k\ell m} e^{k' \ell' m'} = \sum_{abc} g^{k\ell m \ k'\ell'm'}_{\quad abc} 
      e^{abc},
     \qquad
      \Delta(e^{k\ell m}) = \sum_{abc \atop a'b'c'} f_{abc \ a'b'c'}^{\quad  k \ell m} 
      e^{abc} \otimes e^{a' b' c'}.
     \label{AStc} 
\eea
Denoting the generators of $ \A $ as
\[
  e^{100} = x, \qquad e^{010} = z, \qquad e^{001} = y,
\]
we use  (\ref{Uqdef}) and (\ref{Uqcopro}) to read off the structure constants 
in (\ref{UqStc}) and thereby determine the algebraic relations among the 
generators:
\beq
  [x, y] = 0, \qquad [x, z] = 2 \ln q \; x, \qquad [y, z] = 2 \ln q \; y.
  \label{xyzcomm}
\eeq
These are one parameter reduction of the relations in \cite{FG}. 
The full set of the dual basis $ e^{k \ell m} $ may be derived using the 
following recipe. We start by listing the coproduct structure
\bea
\Delta(E_{k\ell m}) &=& \sum_{k' \ell' m'}\bn1{k}{k'}\;\left(\begin{array}{c}
\ell\\ \ell' \end{array}\right)\;\bn1{m}{m'} J_{+}^{k - k'}\,q^{-k' J_{0}}\,
J_{0}^{\ell - \ell'}\,q^{- m' J_{0}}\, J_{-}^{m - m'} \nn\\
& & \otimes J_{+}^{k'}\,q^{(k - k') J_{0}}\,
J_{0}^{\ell'}\,q^{(m- m') J_{0}}\, J_{-}^{m'}, \label{E_co}
\eea    
where the classical and the $q$-binomial coefficients read, respectively:
\beq
\left(\begin{array}{c} n\\ k \end{array}\right) = \frac{n!}{k!\;(n - k)!}, 
\qquad
\bn1{n}{k} = \frac{ [n]_q! }{ [k]_q! [n-k]_q! }. \label{bnc1}
\eeq
Employing (\ref{E_co}) the structure constants $g^{abc \ a'b'c'}_{\ \  k\ell m}$
defined in the second equation in (\ref{UqStc}) may be read explicitly. The 
duality requirement ensures that  these constants compose the algebraic 
multiplication relations of the basis set $e^{k \ell m}$ in (\ref{AStc}). 
Using an inductive procedure we now construct the basis set 
\beq
  e^{k \ell m} = \frac{x^k}{ [k]_q! } 
  \frac{ (z - (k-m) \ln q)^{\ell} }{ \ell! } \frac{ y^m }{ [m]_q! }.
  \label{Abasis}
\eeq
The universal ${\cal T}$-matrix that caps the duality structure of the Hopf 
algebras $\U$ and $\A$ is defined by 
$ \displaystyle {\cal T} = \sum_{k \ell m} e^{k \ell m} \otimes E_{k \ell m},$ 
and its closed form expression may be obtained by substituting the dually 
related basis sets given in (\ref{Ubasis}) and (\ref{Abasis}), respectively:
\beq
 {\cal T} = 
   \left(
     \sum_{k=0}^{\infty} \frac{ (x \otimes J_+ q^{-J_0})^k }{ (k)_{q^{-2}}! }
   \right)
   e^{z \otimes J_0} 
   \left(
     \sum_{m=0}^{\infty} \frac{ (y \otimes q^{J_0} J_- )^m }{ (m)_{q^2} ! }
   \right),
   \label{uT2}   
\eeq 
where
\beq
  (n)_q = \frac{ 1 - q^{n} } {1-q}. \label{qnum2}
\eeq
We emphasize that the classical $ q \rightarrow 1 $ limit of the 
universal ${\cal T}$-matrix (\ref{uT2}) yields the usual exponential 
mapping from the Lie algebra $ su(2) $ to the Lie group $ SU(2). $ 

  The universal ${\cal T}$-matrix (\ref{uT2}) reproduces the standard matrix 
expression of $ \A $ at  the fundamental representation $(\pi)$ of $ \U:$
\beq
  \pi(J_+) = \begin{pmatrix} 0 & 1 \\ 0 & 0 \end{pmatrix},
  \qquad
  \pi(J_-) = \begin{pmatrix} 0 & 0 \\ 1 & 0 \end{pmatrix},
  \qquad
  \pi(J_0) = \begin{pmatrix} \hf & 0 \\ 0 & -\hf \end{pmatrix}.
  \label{UfundR}
\eeq
Denoting the $ 2 \times 2 $ matrix by
\beq
  \begin{pmatrix} a & b \\ c & d \end{pmatrix}
  \equiv
  (id \otimes \pi)({\cal T}) 
  = 
  \begin{pmatrix}
    e^{z/2} + x e^{-z/2} y & q^{1/2} x e^{-z/2} \\
    q^{-1/2} e^{-z/2} y & e^{-z/2}
  \end{pmatrix},
  \label{abcd}
\eeq
one verifies that the elements $ a, b, c, d$ satisfy the familiar defining 
relations of the quantum group $ SL_q(2,{\mathbb C}):$ 
\bea
  & & ab = q^{-1} ba, \qquad ac = q^{-1} ca, \qquad bd = q^{-1} db 
    \nn \\
  & & cd = q^{-1} dc, \qquad bc = cb, \qquad [a, d] = (q^{-1}-q)\, bc, 
    \label{SLq2comm} \\
  & & ad - q^{-1} bc = 1. \nn
\eea
In other words, (\ref{abcd}) is nothing but the Gauss decomposition of 
$ SL_q(2,{\mathbb C}) $ given in \cite{FG}. We remark that $ a, \; d $ are 
invertible, but $ b, \; c$ need not be so. 
Inverting (\ref{abcd}), $ x, y, z $ are expressed in terms of 
$ SL_q(2,{\mathbb C}) $ matrix entries: 
\beq
  x = q^{-1/2} b d^{-1}, \qquad y = q^{1/2} d^{-1} c, \qquad
  e^{z/2} = d^{-1}. 
  \label{xyz-abcd}
\eeq

  The quantum group $ \A $ is a real form of $ SL_q(2,{\mathbb C}). $
As is well-known, the real form $ \A $ is defined by the $*$-involution 
\cite{Wor}
\beq
  a^* = d, \qquad b^* =  - q^{-1} c, \qquad 
  c^* = - q b, \qquad d^* = a.
  \label{SUq2star}
\eeq
The $*$-involution map of the generators of $\A$
obtained via (\ref{xyz-abcd}) and (\ref{SUq2star}) 
\beq 
  x^* = - q^{-1/2} c a^{-1}, \qquad
  e^{z^*/2} = a^{-1}, \qquad y^* = - q^{1/2} a^{-1} b
  \label{xyzstar1}
\eeq
maintains a close kinship to the antipode:
\beq
  S(x) =  y^*, \qquad S(y) =  x^*, \qquad 
  S(z) = z^*. \label{antipode-xyz}
\eeq
Introducing the element
\beq
  \zeta = -q xe^{-z}y = -q bc, \label{zeta}
\eeq
 the $*$-map (\ref{xyzstar1}) is rewritten as 
\beq
  x^* = - \frac{1}{ 1 - \zeta}\, y e^{-z}, 
  \qquad
  e^{z^*/2} = \frac{1}{ 1 - \zeta}\, e^{-z/2},
  \qquad
     y^* = - \frac{1}{ 1 - \zeta}\, e^{-z} x.
     \label{xyzstar2} 
\eeq
One immediately notices that the element $\zeta$ is real: 
$ \zeta^{*} = \zeta. $ 
It is observed that $ x,\; y,\; z $ and their $*$-involutions 
do not commute. For instance, $x$ and $x^*$ satisfy the relations:
\beq
  x x^* = \frac{ q^{-2} x^* x }{ 1 + (q-q^{-1}) x^* x },
  \qquad
  x^* x = \frac{ q^2 x x^* }{ 1 - q^2(q-q^{-1}) x x^* }. 
  \label{xxstar-com}
\eeq
For $j$ being a non-negative integer or a positive half-integer the following 
identities may be proved by induction:
\bea
  & & e^{jz^*} e^{jz} = \frac{1}{ (\zeta;q^2)_{2j} }, 
     \qquad \qquad \quad
     e^{jz}\, e^{jz^*} = \frac{1}{ (q^{-2} \zeta;q^{-2})_{2j} }, 
    \label{ez1} \\[5pt]
  & & e^{-jz^*} e^{-jz} =  (q^{-2} \zeta;q^{-2})_{2j},
    \qquad
    e^{-jz}\, e^{-jz^*} = (\zeta;q^2)_{2j}, \label{ez2}
\eea
where $ (a;q)_n $ is the $q$-shifted factorial. For a positive integer $n$
the identities given below hold:
\beq
  (x^*)^n x^n = \frac{ q^{n(n-2)} \zeta^n }{ (\zeta;q^2)_n },
  \qquad
  x^n (x^*)^n = \frac{  q^{-n(n+2)} \zeta^n}{ (q^{-2} \zeta; q^{-2})_n }.
  \label{xnxsn}
\eeq
We will use these relations subsequently. Definitions and formulae of 
$q$-analysis used in this work are summarized in Appendix.

 The coproduct of $\A$ can be regarded as a left/right coaction 
of $ \A $ on $ \A $ itself. This leads a left and a right coaction of 
$ \A $ on $ x $ and $ y, $ respectively.  Denoting these coactions by  
$ \varphi(x),\;\varphi(y) $ we have
\bea
 \varphi(x) &=& q^{-1/2} \Delta(bd^{-1}) \nn \\
    &=& (a \otimes x + b \otimes q^{-1/2}) (q^{1/2} c \otimes x + d \otimes 1)^{-1},
    \label{Lactionx} \\[3pt]
 \varphi(y) &=& q^{1/2} \Delta(d^{-1} c) \nn \\
    &=& (q^{-1/2}y \otimes b + 1 \otimes d)^{-1} (y \otimes a + q^{1/2} \otimes c).
    \label{Ractiony}
\eea
It is easy to verify that 
(\ref{xxstar-com}) is covariant under the coaction (\ref{Lactionx}).

%%%%%%%%%%%%%%%%%%%%%%%%%%%%%%%%%%%%%%%%%%%%%%%%%%%%%%%%%
%
%  SU_q(2) coherent states
%
%%%%%%%%%%%%%%%%%%%%%%%%%%%%%%%%%%%%%%%%%%%%%%%%%%%%%%%%%
%
\section{\bm{$SU_q(2)$} coherent states}
\setcounter{equation}{0}
\label{SUq2CS}

\subsection{Definition of coherent states}
\label{subsecDef}

  The universal ${\cal T}$-matrix (\ref{uT2}) was constructed 
as a quantum analogue of an exponential mapping. 
As in \cite{JS},  the coherent state for $\A$ is defined as a quantum analogue 
of the generalized coherent state. Namely, we consider 
a representation of $ \U $ and take a fiducial vector from the representation space. 
The coherent state 
is the state obtained by transformation of the fiducial vector by 
the universal ${\cal T}$-matrix. 
Let us take spin $j$ representation of $ \U:$
\bea
  & &  J_{\pm} \ket{jm} = \sqrt{ [j \mp m]_q[j \pm m + 1]_q } \ket{j\; m\pm 1}, 
    \nn \\[5pt]
  & & J_0  \ket{jm} = m \ket{jm}, 
    \label{SUq2rep} 
\eea
where $j$ is a non-negative integer or a positive half-integer. 
This is a unitary representation for $ q \in U(1) $ or $ q \in {\mathbb R}. $ 
We assume $ q \in {\mathbb R} $ throughout this article. 
Coherent state for $ \A $ is defined by
\beq
  \ket{x, z} = {\cal T} \ket{j \; -j}.    \label{CSdef}
\eeq
Repeated use of (\ref{SUq2rep}) gives us the explicit form of the coherent state:
\bea
   \ket{x, z} &=& 
    \sum_{n=0}^{2j} q^{nj} \bn1{2j}{n}^{1/2} x^n e^{-jz} \ket{j\; -j+n}
     \label{CS2-1} \\
   &=&
   e^{-jz} \sum_{n=0}^{2j} q^{-nj} \bn1{2j}{n}^{1/2} x^n \ket{j\; -j+n}.
     \label{CS2-2}
\eea
In the classical case if two coherent states differ from one another only by a 
phase factor, they are regarded as the same state. Thus the coherent state has 
one-to-one correspondence to a point of the coset space $ SU(2)/U(1). $ 
Standard choice \cite{Per,ACGT,CEP} of the coset representative with $z$ being 
real, is given by 
\beq
  e^z = 1 + |x|^2, \qquad 
  x^* = -y, \qquad y^* = - x. \label{z-real}
\eeq
In the quantum case, however, a choice of real $ z $ is inadmissible since, as 
may be seen from (\ref{xyzstar2}), it leads to a contradiction: 
\[
  (1- q^{-2})\, \zeta = 0.
\]
We thus take the universal $ \cal T$-matrix itself 
as our coset representative. 

  By definition, the coherent state (\ref{CSdef}) has unit norm. 
Proof by direct computation is also easy. We give the proof below 
in order to see that the factor $ e^{-jz} $ behaves classically. 
With the aid of the identities (\ref{qid2}) and (\ref{bn-qfac}) for 
$q$-shifted factorials, we evaluate the norm:
\[
   \langle x, z | x, z \rangle
   = e^{-jz^*} e^{-jz} \,
       \sum_{n=0}^{2j} (-1)^n q^{n(n-1)} 
       \frac{ (q^{-4j};q^2)_n }{ (q^2;q^2)_n (q^{-4j}\zeta;q^2)_n }\, \zeta^n.
\]
It is remarkable that $ e^{-jz^*} e^{-jz} $ is factored out and, because of 
(\ref{ez2}), commutes with the rest of the expression. This precisely happens 
in the classical case. Expressing the norm in terms of the basic hypergeometric
function and using (\ref{11trans}), it is observed to be equal to unity:
\bea
  \langle x, z | x, z \rangle
  &=& e^{-jz^*} e^{-jz} \,
     \BH{1}{1}\BHA{q^{-4j}}{q^{-4j}\zeta}{q^2}{\zeta} 
     \nn \\[5pt]
  &=& e^{-jz^*} e^{-jz} \, 
      \frac{ (\zeta;q^2)_{\infty} }{ (q^{-4j} \zeta;q^2)_{\infty} } = 1.
     \nn
\eea

%%%%%%%%%%%%%%%%%%%%%%%%%%%%%%%%%%%%%%%%%%%%%%%%%%%%%%%%%
%
\subsection{Resolution of unity}
\label{subsecRU}

  One of the most important properties of the generalized coherent 
states is "resolution of unity". This is essential for applications 
of coherent states to physical problems, path integrals, representation 
theory, and so on. Despite the noncommutativity of $ \A, $ 
the coherent state (\ref{CS2-1}) satisfies the resolution of 
unity with respect to an invariant integration over quantum groups. 
Invariant integration over the quantum group $ \A $ has been discussed in 
\cite{Wor,Woro,MMNM,APVV}. Attempt was made in \cite{APVV} to develop a 
semiclassical approach towards the construction of the quantum Haar measure.
In our work we follow the description given in \cite{CP}.

  Let ${\cal G}$ be an arbitrary quantum group. A linear functional 
$ H : {\cal G} \ \longrightarrow \ {\mathbb C} $ is said to be 
\textit{normalized bi-invariant integral} if 
\begin{enumerate}
  \item $ \Haar{1_{\cal G}} = 1,$
  \item for any $ f \in {\cal G} $
  \beq
      (H \otimes id)[\; \Delta(f) \; ] = (id \otimes H)[ \; \Delta(f)\; ]
        = \Haar{f}. \label{Hdef}
  \eeq 
\end{enumerate}
Let ${\cal V}$ be an algebra dual to $ {\cal G}. $  
Writing the coproduct of $ f \in {\cal G} $ as 
$ \displaystyle \Delta(f) = \sum_k f_k \otimes f^k, $ 
the left and the right action of $ Z \in {\cal V} $ on $f$ are 
defined by
\beq
  Z \rhd f = \sum_k f_k \langle Z, f^k \rangle, 
  \qquad
  f \lhd Z = \sum_k \langle Z, f_k \rangle f^k.
  \label{LRaction-def}
\eeq
Then one can prove the left and the right invariance of the integral $H:$
\beq
  \Haar{ Z \rhd f } = H[\; f \lhd Z \;] = \epsilon(Z)\, \Haar{f}. 
  \label{LRinvofH}
\eeq
Below we cite some results for $ {\cal G} = SL_q(2,{\mathbb C}), \ {\cal V} 
= U_q[sl(2,{\mathbb C})] $  from \cite{CP}, since they are used in the proof 
of resolution of unity (note that our 
conventions are slightly different from the ones in \cite{CP}). 
From the invariance under the action of $ q^{J_0},$ it is shown that
\bea
  & & 
  \Haar{ a^{\alpha} b^{\beta} c^{\gamma} d^{\delta} } \neq 0, 
  \quad \mbox{only if} \quad  \alpha = \delta \ \mbox{and} \ \beta = \gamma,
  \label{HSL2-1} \\[5pt]
  & & 
  \Haar{ b^{\beta} c^{\gamma} d^{-\delta} a^{-\alpha} } \neq 0, 
  \quad \mbox{only if} \quad  \alpha = \delta \ \mbox{and} \ \beta = \gamma.
  \label{HSL2-2}
\eea
Implementing unit value of the quantum determinant the monomials of the form
$ a^{\alpha} b^{\beta} c^{\beta} d^{\alpha} $ and  
$ b^{\beta} c^{\beta} d^{-\alpha} a^{-\alpha} $ are 
reduced to a power series of the element $ \zeta $ defined in (\ref{zeta}). 
Integration of $ \zeta^n $ is computed using bi-invariance 
(\ref{Hdef}):
\beq
  \Haar{\zeta^n} = \frac{q^{2n}}{ (n+1)_{q^2}}.
  \label{Hzn}
\eeq
It is worth noting that that for a polynomial $ f(\zeta) $ 
the invariant integral can be written as a  $q$-integral:
\beq
  \Haar{ f(\zeta) } = \int_0^1 f(q^2 \zeta)\, d_{q^2}\zeta. \label{Hbyint}
\eeq

 With these preparations, one can show that the coherent state for quantum 
group $ \A$ provides the resolution of unity:  
\beq
  (2j+1)_{q^2} H\big[\; \ket{x,z} \bra{x,z} \;\big] = 1. 
  \label{SU2unity}
\eeq 

\noindent
{\bf Proof:} By virtue of (\ref{ez2}), the integrand is written as:
\bea
   \ket{x,z} \bra{x,z}&=&
     \sum_{n,m=0}^{2j} q^{j(n+m)} \bn1{2j}{n}^{1/2} \bn1{2j}{m}^{1/2}
    \nn \\[4pt]
  & & \qquad \times \;
     (q^{-2n}\zeta;q^2)_{2j} \, x^n (x^*)^m \, \ket{j\; -j+n} \bra{j \; -j+m}. 
    \label{integrand}
\eea 
Noting that 
$ x^n (x^*)^m = (-1)^m q^{(m^2-n^2)/2-m-nm}\, b^n c^m d^{-n} a^{-m}, $
we see with the aid of (\ref{HSL2-2}) that the terms in (\ref{integrand}) 
contribute to the integration (\ref{SU2unity}) only if $ n=m:$ 
\bea
  & & \hspace{-2cm} H\big[\; \ket{x, z} \bra{x, z} \;\big] 
    \nn \\
  & \quad = \ &
   \sum_{n=0}^{2j} q^{2jn} \bn1{2j}{n} 
      H[\; (q^{-2n}\zeta;q^2)_{2j} \, x^n (x^*)^n \;] 
      \ket{j\; -j+n} \bra{j \; -j+n}
  \nn \\
  &  \quad \stackrel{(\ref{xnxsn})}{=} \ &
   \sum_{n=0}^{2j} q^{2jn -n(n+2)} \bn1{2j}{n} 
      H[\; (\zeta;q^2)_{2j-n} \zeta^n \;] 
      \ket{j\; -j+n} \bra{j \; -j+n}.
   \nn
\eea
Before proceeding further, we need to show two identities:
\bea
  & & (\zeta;q)_m = 
    \sum_{k=0}^m q^{mk} \, \frac{ (q^{-m};q)_k }{ (q;q)_k } \, \zeta^k, 
    \label{qrel1} \\[5pt]
  & & \sum_{k=0}^{m-n} \frac{ q^{(m-n+1)k} }{ (n+k+1)_q } 
      \frac{ (q^{-m+n};q)_k }{ (q;q)_k } 
      =
      (-1)^n \, \frac{ q^{-mn + n(n-1)/2} }{ (m+1)_q } 
      \frac{ (q;q)_n }{ (q^{-m};q)_n }.
      \label{qrel2}
\eea
The identity (\ref{qrel1}) is obtained by setting $ a = q^{-m}, 
\ z = q^m \zeta $ in the $q$-binomial theorem (\ref{q-bThm}). 
To prove (\ref{qrel2}), we start with a slight modification of (\ref{qrel1}):
\[
   \zeta^n (\zeta;q)_{m-n} = 
    \sum_{k=0}^{m-n} q^{(m-n)k} \, \frac{ (q^{-m+n};q)_k }{ (q;q)_k } \, \zeta^{n+k},
\]
and perform the $q$-integration on both sides. The rhs immediately yields
\[
  I_{m-n,n} \equiv \int_0^x 
    \sum_{k=0}^{m-n} q^{(m-n)k} \, \frac{ (q^{-m+n};q)_k }{ (q;q)_k } \, \zeta^{n+k} 
  d_q \zeta 
  = 
  \sum_{k=0}^{m-n} q^{(m-n)k} \, \frac{ (q^{-m+n};q)_k }{ (q;q)_k } \,
  \frac{ x^{n+k+1} }{ (n+k+1)_q }.
\]
We use integral by parts (\ref{intpart}) for the lhs:
\bea
  I_{m-n,n} &=&
  \frac{ x^{n+1} }{ (n+1)_q } (xq^{-1};q)_{m-n} - 
  \frac{1}{ (n+1)_q } \int_0^x \zeta^{n+1} D_q (\zeta q^{-1};q)_{m-n} d_q \zeta
  \nn \\[5pt]
  &=& 
  \frac{ x^{n+1} }{ (n+1)_q } (xq^{-1};q)_{m-n} + 
  \frac{ q^{-1} (m-n)_q }{ (n+1)_q } \, I_{m-n-1,n+1}.
  \nn  
\eea
As the first term vanishes for the intended value $x=q$, we do not keep track 
of it. Consequently the recurrence relation is easily solved: 
\bea
  I_{m-n,n} &=& q^{-(m-n)} 
    \frac{ (m-n)_q! }{ (n+1)_q (n+2)_q \cdots (m)_q } \, I_{0,m}
   \nn \\[5pt]
  &=&
    (-1)^n q^{-(m-n)-mn + n(n-1)/2} 
    \frac{ (q;q)_n }{ (q^{-m};q)_n } 
    \frac{ x^{m+1} }{ (m+1)_q }.
  \nn
\eea
Setting $ x = q, $ we obtain (\ref{qrel2}). 

Returning to the resolution of unity, one can now  
carry out the integration:
\bea
  & & H\big[\; \ket{x,z} \bra{x,z} \;\big] 
    \nn \\
  & & \quad \stackrel{(\ref{qrel1})}{=} \ 
   \sum_{n=0}^{2j} (-1)^n q^{4jn -n(n+1)} \,
   \frac{ (q^{-4j};q^2)_n }{ (q^2;q^2)_n } 
   \sum_{k=0}^{2j-n} q^{2(2j-n)k} \,
   \frac{ (q^{-4j+2n};q^2)_k }{ (q^2;q^2)_k }  
   \nn \\[5pt]
  & & \quad \;  \times \; 
      H[\; \zeta^{n+k} \;] 
      \ket{j\; -j+n} \bra{j \; -j+n}
   \nn \\[5pt]
  & & \quad \stackrel{(\ref{Hzn})}{=} \ 
   \sum_{n=0}^{2j} (-1)^n q^{4jn -n(n-1)} \,
   \frac{ (q^{-4j};q^2)_n }{ (q^2;q^2)_n } 
   \sum_{k=0}^{2j-n} \frac{ q^{2(2j-n+1)k} }{ (n+k+1)_{q^2} } 
   \frac{ (q^{-4j+2n};q^2)_k }{ (q^2;q^2)_k }
   \nn \\
  & & \quad \; \times \;
      \ket{j\; -j+n} \bra{j \; -j+n}
   \nn \\[5pt]
  & & \quad \stackrel{(\ref{qrel2})}{=} \ 
    \frac{1}{ (2j+1)_{q^2} } \sum_{n=0}^{2j}\, \ket{j\; -j+n} \bra{j \; -j+n} 
    =  \frac{1}{ (2j+1)_{q^2} }.
  \nn
\eea
We thus proved (\ref{SU2unity}). $\Box$

%%%%%%%%%%%%%%%%%%%%%%%%%%%%%%%%%%%%%%%%%%%%%%%%%%%%%%%%
%
\subsection{Properties of coherent state}
\label{subsecProp}

The coherent state for the quantum group $ \A $ allows for easy generalizations
of some properties whose classical counterparts are well-known \cite{Per}. 
The discussions presented here in conjunction with the results obtained in the 
preceding and the succeeding sections show that many key characteristics of the
$ SU(2) $ coherent state can be lifted up to the quantum group setting. 

  We first investigate the overlap of two coherent states. The coherent state 
$ \ket{x,z} $ is an element in $ \A \otimes V^{(j)}, $ where $ V^{(j)} $ is 
the representation space of spin $j.$ According to \cite{JS}, 
the overlap of two coherent states is defined as an object in $\A \otimes \A.$ 
Let us introduce two independent copies of the coherent state distinguished by 
the subscripts: $ \ket{x_a, z_a}, \ a = 1, 2. $ Then the recipe given in
\cite{JS} yields the overlap as
\beq
  \bracket{x_1, z_1}{x_2, z_2} = 
  ( e^{-jz_1^*} \otimes e^{-jz_2} ) \sum_{n=0}^{2j} 
  \bn1{2j}{n} (x_1^* \otimes x_2)^n.
  \label{overlap}
\eeq
In (\ref{overlap}) the exponential terms contained in the parenthesis commute
with the factored sum. 

  Next we enlist the actions of the generators of the $ \U $ algebra on the 
coherent state:
\bea
  & & J_+ \ket{x,z} = \sum_{n=0}^{2j} q^{(n-1)j} 
      \qn{n} \bn1{2j}{n}^{1/2}  x^{n-1} e^{-jz} \ket{j\; -j+n},
    \label{JponCS} \\[3pt]
  & & J_- \ket{x,z} = \sum_{n=0}^{2j} q^{(n+1)j} 
      \qn{2j-n} \bn1{2j}{n}^{1/2}  x^{n+1} e^{-jz} \ket{j\; -j+n},
    \label{JmonCS} \\[3pt]
  & & \qn{J_0} \ket{x,z} = \sum_{n=0}^{2j} q^{nj} 
      \qn{-j+n} \bn1{2j}{n}^{1/2}  x^{n} e^{-jz} \ket{j\; -j+n}.
    \label{J0onCS}
\eea
It now follows that there exists an operator which annihilates 
the coherent state:
\beq
  ( J_- + (1+q^{2j})x \qn{J_0} - q^{2j} x^2 J_+ ) \ket{x,z} = 0.
  \label{CSanni}
\eeq
Furthermore, the coherent state is an eigenvector of the 
operator $ \Gamma $ defined by
\beq
  \Gamma = 
  (1 - (q^j + q^{-j}) \zeta )\, q^{J_0} \qn{J_0} 
   - (1 - q^j \zeta)\, x \, q^{J_0} J_+ + q^{-j-1} e^{-z} y\, q^{J_0} J_-.
  \label{Gammadef}
\eeq
Eigenvalue relation reads
\beq
  \Gamma \ket{x,z} = - q^{-j} \qn{j} \ket{x,z}. 
  \label{GonCS}
\eeq
Relations (\ref{CSanni}) and (\ref{GonCS}) may be proven by 
straightforward computation.

%%%%%%%%%%%%%%%%%%%%%%%%%%%%%%%%%%%%%%%%%%%%%%%%%%%%%%%%%
%
%  Coherent state and q-sphere
%
%%%%%%%%%%%%%%%%%%%%%%%%%%%%%%%%%%%%%%%%%%%%%%%%%%%%%%%%%
%
\section{Coherent state and $q$-sphere}
\setcounter{equation}{0}
\label{CSqsphere}

  Geometrical importance of coherent states is due to 
the fact that they provide natural description 
of K\"ahler structure of homogeneous spaces \cite{Ono}. 
Let us recall the case of $ SU(2). $ The homogeneous space 
for the $ SU(2) $ coherent state is  $ SU(2)/U(1) \approx S^2. $ 
Expectation values of the $su(2)$ Lie algebra elements with respect to 
the $ SU(2)$ coherent state reflect this fact:
\beq
  \langle J_+ \rangle \langle J_- \rangle + \langle J_0 \rangle^2 
  = 
  \langle J_x \rangle^2 + \langle J_y \rangle^2 
  + \langle J_z \rangle^2 = j^2.
  \label{S2class}
\eeq
Since the expectation values $ \langle J_a \rangle  $ are functions on a  
complex plane, this provides a complex description of the 2-sphere: 
$ S^2 \approx {\mathbb C} \cup \{ \ \infty \ \}. $ 
K\"ahler potential on $ S^2 $ is given by the normalization factor 
of the coherent state: 
$ F(x,x^*) = -\ln |\bracket{j-j}{x}|^2. $ 

  Now we turn to the quantum group setting. 
Homogeneous spaces for quantum groups are introduced in \cite{DK}. 
Making use of the quantum subgroups of $ \A $ consisting of 
the set of diagonal matrices 
$ U(1) = \{ \ diag(\alpha,\alpha^*) \ | \ |\alpha|=1, \ \alpha 
\in {\mathbb C} \}, $ 
one can see that the homogeneous space $ SU_q(2)/U(1) $ is 
generated by $ ab,\; bc,\; cd. $ Thus the homogeneous space $ SU_q(2)/U(1) $ 
is identified with Podle\'s $q$-sphere \cite{Pod} embedded into 
$ \A $ \cite{DK,Koo}:
\beq
  x_{-1}  = \sqrt{1+q^2}  \,  ab, \qquad 
  x_0 = 1 + (q + q^{-1}) \,  bc \qquad 
  x_1 = \sqrt{1+q^{-2}}  \, dc. \label{qsphere-SL/K}
\eeq
These coordinates of $q$-sphere  satisfy the relations
\bea
  & & x_0^2 - q^{-1} x_1 x_{-1} - q x_{-1} x_1 = 1, \nn \\[3pt]
  & & (1-q^{-2})  \, x_0^2 + q^{-1} x_{-1} x_1 - q^{-1} x_1 x_{-1} = (1-q^{-2})  \, x_0, 
     \label{qsphere-def} \\[3pt]
  & & x_{-1} x_0 - q^{-2} x_0 x_{-1} = (1-q^{-2})  \, x_{-1}, \nn \\[3pt]
  & & x_0 x_1 - q^{-2} x_{1} x_0 = (1-q^{-2})  \, x_1. \nn
\eea
This is a special case of the $q$-sphere $ S^2_{q,\rho} $ 
for a specific value of the  parameter $ \rho.$   
Infinitesimal characterization \cite{Koo,DK} for this $q$-sphere is given by
\beq
  (q^{J_0} - q^{-J_0}) \rhd x_k = 0, \qquad k = \pm 1, 0. \label{inf-char}
\eeq

  Let us consider expectation values of some specific elements 
in $ \U $ with respect to the coherent state for $ \A: $ 
\bea
  & & X_+ \equiv \bra{x,z} J_+ q^{-J_0} \ket{x,z}
      = \qn{2j} (1 - \zeta) \, x^{*},
   \nn \\[3pt]
  & & X_- \equiv \bra{x,z} q^{-J_0} J_-  \ket{x,z}
      = \qn{2j} x \, (1 - \zeta),
   \label{Xm} \\[3pt]
  & & X_0 \equiv \bra{x,z} q^{-J_0} \qn{J_0} \ket{x,z}
      = q^{-2} \qn{2j} \zeta - q^j \qn{j}.
   \nn
\eea
These expectation values, after suitable scaling given by 
\bea
  & & x_1 = -q \frac{ \sqrt{\qn{2}} }{ \qn{2j} } X_+ 
          = - q \sqrt{ \qn{2} }  (1-\zeta) x^*,
   \nn \\[3pt]
  & & x_0 = 1 - q \frac{ \qn{2} }{ \qn{2j} }(X_0 + q^j\qn{j})
          = 1 - q^{-1} \qn{2} \zeta,
   \label{x-Xrelation} \\[3pt]
  & & x_{-1} = \frac{ \sqrt{\qn{2}} }{ \qn{2j} } X_-
             = \sqrt{ \qn{2} } x (1-\zeta),
   \nn 
\eea
satisfy the defining relations (\ref{qsphere-def}) of the $q$-sphere,
while maintaining the following $*$-involution map:
\beq
  x_1^* = -q x_{-1}, \qquad x_0^* = x_0, \qquad x_{-1}^* = -q^{-1} x_1.
  \label{star-sphere}
\eeq
The relation
\beq
 \zeta = \frac{ q x^* x }{ 1 + q x^* x } = 
         \frac{ q^3 x x^* }{ 1 + q x x^* }
 \label{zeta-xx}
\eeq
ensures that the coordinate $ x_k $ are functions of only $ x $ and $ x^*. $ 
Thus the coherent state gives a natural complex description of the
$q$-sphere. We emphasize that we did not introduce any additional assumption 
to obtain the above complexification. It is a direct consequence of the 
definitions of $ \U $ and its dual together with the finite dimensional 
representations of $ \U. $  Complex description of the quantum 2-sphere has 
been previously considered in \cite{St,CHZ}. In \cite{St}, the dressing 
transformation \cite{JS2} is used to 
derive the relations between the complex coordinates $ x,\; x^* $ of 
$q$-sphere. It looks a bit more complex than our relation 
(\ref{xxstar-com}). The authors of \cite{CHZ} introduced $ b^{-1},\; c^{-1} $ 
to derive the stereographic projection of the $q$-sphere. The resulting 
commutation relations of complex coordinates are rather simple. 
However, $ b,\; c \in \A $ are not assumed to be invertible in our setting.

  Next we study differential calculus on 
$q$-sphere in our complex description. Recall that differential calculus 
on $q$-sphere in $ x_{\pm 1}, \; x_0 $ coordinates have been 
extensively studied \cite{Pod2,Pod3,AS}. Our aim in this section is to 
develop differential calculus in the complex coordinates $x,\; x^*. $ Such a 
differential calculus is considered in \cite{CHZ} based on different 
noncommutative coordinates from ours. 

  As seen in (\ref{qsphere-SL/K}), the $q$-sphere can be embedded in $ \A. $ 
The embedding allows us to infer the differential structures on the $q$-sphere 
from the well-known \cite{Wor,Wor2} covariant differential calculi on $ \A. $ 
We use the so-called left-covariant 3D calculus \cite{Wor} on 
$ \A $ because of the reasons advanced in \cite{CHZ}. 
Along the line in \cite{SS}, we list the relations of the  
3D calculus on $ \A $ in our conventions for subsequent use. 
Three following elements in $ \U$
\beq
  {\cal X}_1 = \frac{ q^{4J_0} - 1 }{ q^2 - 1 }, \qquad
  {\cal X}_0 = q^{1/2} J_+ q^{J_0}, \qquad
  {\cal X}_2 = q^{-1/2} J_- q^{J_0}. 
  \label{QtanSp}
\eeq
may be regarded to span the quantum tangent space on $ \A,$ and let
$ \omega_k \ (k = 0, 1, 2) $ be the one-forms dual to these tangent vectors,
respectively. The differentials of the elements $ a, b, c, d, $ denoted  
sequentially as $ \alpha, \beta, \gamma, \delta, $
relate to $ \omega_k $ as follows:
\beq
  \alpha = a \omega_1 + b \omega_2, \quad 
  \beta = a \omega_0 - q^{-2} b \omega_1, \quad
  \gamma = c \omega_1 + d \omega_2, \quad 
  \delta = c \omega_0 - q^{-2} d \omega_1.
\eeq
Commutation relations between the coordinates and differentials read:
\bea
  & & \omega_1 a = q^2 a \omega_1, \qquad \omega_1 b = q^{-2} b \omega_1, 
      \qquad
      \omega_1 c = q^2 c \omega_1, \qquad \omega_1 d = q^{-2} d \omega_1,
  \nn \\[3pt]
  & & \omega_k a = q a \omega_k, \qquad \omega_k b = q^{-1} b \omega_k, 
      \qquad 
      \omega_k c = q c \omega_k, \qquad \omega_k d = q^{-1} d \omega_k,
  \label{3Dcom}
\eea
where $ k = 0, 2. $ It follows that
\beq
  [x,\; \omega_k ] = [x^*, \; \omega_k] = 0, 
  \qquad k = 0, 1, 2.
  \label{x-omega-rel}
\eeq

  With these settings, it is an easy exercise to show that
\bea
 & & x dx = q^2 dx \, x, \qquad x^* dx^* = q^{-2} dx^*\, x^*, 
   \nn \\[5pt]
 & & dx \, x^* = q^{-2} f_-(x^*x) x^* dx, \quad
     dx^* \, x = q^2 x f_+(x^*x) dx^*,
   \label{x-dx}
\eea
where
\beq
  f_{\pm}(x^*x) = \frac{1 - \zeta}{1-q^{\pm 4} \zeta}
  = \frac{1}{1 + (1-q^{\pm 4})\, q\, x^* x}. \label{fpm}
\eeq
The nilpotency of the complex differentials follow:
\beq
   (dx)^2  = (dx^*)^2 = 0. \label{nilpo-dx}
\eeq
To determine the commutation relation between $ dx $ and $ dx^*, $ 
we need to calculate $ df_{\pm}(x^* x). $ 
This is done with the help of identities:
\bea
 & & dx \, t = q^{-4}\, t f_-(t) \, dx, \qquad
     dx^* \, t = q^4 \, t f_+(t) \, dx^*, \nn \\[3pt]
 & & x \, t f_+(t) = \frac{ q^{-2} t }{ 1 - q^2 \omega t } \, x, \qquad
     x^* \, t f_-(t) = \frac{ q^2 t }{ 1 + \omega t } \, x^*, \nn 
\eea
where $ t = x^* x $ and $ \omega = q - q^{-1}. $ By induction, one can 
show that
\beq
  d ( x^* x )^n = \frac{q}{\omega} \frac{1 + \omega t}{1+ qt} \,
       t^{n-1} 
       \left\{
         - \left(
             1 - \frac{q^{2n}}{(1-q^2 \omega t)^n}  
           \right) x\, dx^*
         + \left(
             1 - \frac{q^{-2n}}{(1+\omega t)^n}
           \right) x^* \, dx
       \right\}.
  \label{dxn}
\eeq
After a lengthy computation, we derive the relation:
\beq
  dx \, dx^* = - \frac{f_-(x^*x)}{q^2} 
  \frac{ 1 - q^2 \omega x^*x }{ 1+ \omega x^* x } 
  \, dx^*\, dx.
  \label{dxdxstar}
\eeq
This completes our derivation of the differential calculus on $q$-sphere  
subject to the complexification adopted here.

%%%%%%%%%%%%%%%%%%%%%%%%%%%%%%%%%%%%%%%%%%%%%%%%%%%%%%%%%
%
%  Coherent state representation of U_q[su(2)]
%
%%%%%%%%%%%%%%%%%%%%%%%%%%%%%%%%%%%%%%%%%%%%%%%%%%%%%%%%%
%
\section{Coherent state representation of \bm{$U_q[su(2)]$}}
\setcounter{equation}{0}
\label{CSrepU}

  As an application of the resolution of unity (\ref{SU2unity}),
we discuss a representation of the algebra $ \U $ by making use of 
the coherent state (\ref{CS2-1}). Let $ \ket{c} $ be an arbitrary state 
in the representation space of spin $j:$
\[
  \ket{c} = \sum_{m=-j}^j c_m \ket{jm}, \qquad c_m \in {\mathbb C}
\]
Then by virtue of the resolution of unity, it follows
\bea
  \ket{c} &=& 
  \sum_{m=-j}^j c_m\, (2j+1)_{q^2} \Haar{ \ket{x,z} \bra{x,z} } \ket{jm}
  \nn \\
  &=& (2j+1)_{q^2} 
   \Haar{ \ket{x,z} \sum_{m=-j}^j c_m \bracket{x,z}{jm} },
  \nn
\eea
where
\beq
  \bracket{x,z}{jm} =  
  q^{j(j+m)} \bn1{2j}{j+m}^{1/2} e^{-jz^*} (x^*)^{j+m}.
  \nn %\label{exp-coeff}
\eeq
This implies that any state in spin $j$ representation of $ U_q[su(2)] $ 
can be expanded in the coherent state basis. 
The expansion coefficient is a 
polynomial in $ x^* $ with degree up to $2j. $ 

 Let us consider the monomials
\beq
  \Psi^j_m(x) = q^{-j(j+m)-(j-m)} \bn1{2j}{j+m}^{1/2} x^{j+m},
    \label{basisH}
\eeq
where $ m = -j, -j+1, \cdots, j. $ 
We shall show that these monomials span a vector space carrying a unitary 
representation of $ \U $ with spin $j.$ For arbitrary elements 
$ f(x), \; g(x) $ of this vector space of the monomials we adopt
the following definition of an inner product:
\beq
  ( f(x), g(x) ) 
  = (2j+1)_{q^2} \Haar{ f(x)^* e^{-jz^*} e^{-jz} g(x) }.
  \label{innerH}
\eeq
The monomials (\ref{basisH}) form an orthonormal basis with respect to 
the inner product:
\beq
  ( \Psi^j_{m'}(x), \Psi^j_m(x) ) = \delta_{m'm}. \label{normalH}
\eeq
The proof of the relation (\ref{normalH}) is summarized below. The requirement
(\ref{HSL2-2}) ensures that the rhs of (\ref{normalH}) vanishes 
for $ m' \neq m. $ 
Consequently, we just compute the case of $ m' = m. $ 
With the aid of $q$-binomial theorem (\ref{q-bThm}), 
one verifies that
\bea
  (x^*)^{j+m} e^{-jz^*} e^{-jz} x^{j+m} &=& 
  q^{(j+m) (j+m-2)} \zeta^{j+m} \frac{ (q^{-2(j-m)}\zeta;q^2)_{\infty} }{ (\zeta;q^2)_{\infty} } 
  \nn \\[3pt]
  &=& 
  q^{(j+m) (j+m-2)} \sum_{k=0}^{j-m} \frac{ (q^{-2(j-m)};q^2)_{k} }{ (q^2;q^2)_{k} } \zeta^{j+m+k}. 
  \nn
\eea
Then the invariant integration is carried out easily
\bea
  & & \Haar{ (x^*)^{j+m} e^{-jz^*} e^{-jz} x^{j+m} } 
  \nn \\[3pt]
  & & \qquad   = 
  q^{(j+m) (j+m-2)} \sum_{k=0}^{j-m} \frac{ (q^{-2(j-m)};q^2)_{k} }
  { (q^2;q^2)_{k} } 
  \frac{ q^{2(j+m+k)} }{ (j+m+k+1)_{q^2} }
  \nn \\[3pt]
  & & \qquad = 
  q^{2j(j+1) + 2m(j-1)} \bn1{2j}{j+m}^{-1} \frac{1}{ (2j+1)_{q^2} },
  \label{innerHpr1}
\eea
where the last equality is due to the identity
\beq
  \sum_{k=0}^{j-m} \frac{ (q^{-(j-m)};q)_k }{ (q;q)_k } 
  \frac{ q^{j+m+k} }{ (j+m+k+1)_q } 
  = 
  q^{j^2-m^2} \, \frac{ (j+m)_q ! (j-m)_q ! }{ (2j)_q ! } 
  \frac{ q^{2j} }{ (2j+1)_q }.
  \label{identityforH}
\eeq
Using (\ref{intpart2}) and (\ref{der-form2}), the identity (\ref{identityforH})
is proved in a way parallel to (\ref{qrel2}). Then the normalization condition 
$ ( \Psi^j_{m}(x), \Psi^j_m(x) ) = 1 $  
follows immediately from (\ref{innerHpr1}). 

  Now we introduce a formal derivative operator $ \partial $ by 
\beq
  \partial x^n = n x^{n-1}, \label{formalder}
\eeq
and employ the $q$-derivative to realize the quantum algebra $ \U $ on the 
space spanned by $ \Psi_m^j(x). $ Towards this end we define the operators
\bea
  & & J_+ = ( (2j)_{q^2} x  - x^2 D_{q^2} ) \, q^{3/2-3j-J_0}, 
      \qquad
      J_- = D_{q^2} q^{-1/2+j-J_0},
  \label{suq2onpoly} \\[5pt]
  & & J_0 = x \partial - j. \nn
\eea
Their actions on the monomial basis set read
\bea
  & & J_0 \Psi^j_m = m \Psi^j_m, \nn \\[5pt]
  & & J_{\pm} \Psi^j_m = \sqrt{ (j \mp m)_{q^2} (j \pm m + 1)_{q^2} } \; \Psi^j_{m \pm 1}.
    \label{Jonpoly}
\eea
The hermiticity of the representation (\ref{Jonpoly}) is ensured  by the 
relations 
\[
  (J_0^* \Psi^j_{m'}, \Psi^j_{m}) = (J_0 \Psi^j_{m'}, \Psi^j_{m}),
  \qquad
  (J_{\pm}^* \Psi^j_{m'}, \Psi^j_{m}) = ( J_{\mp} \Psi^j_{m'}, \Psi^j_{m}),
\]
which may be easily proven by straightforward computation.

%%%%%%%%%%%%%%%%%%%%%%%%%%%%%%%%%%%%%%%%%%%%%%%%%%%%%%%%%
%
%  Contraction and entanglement
%
%%%%%%%%%%%%%%%%%%%%%%%%%%%%%%%%%%%%%%%%%%%%%%%%%%%%%%%%%
%
\section{High spin limit}
\setcounter{equation}{0}
\label{Contra-and-Ent}

  In this section, we 
%investigate two more properties 
%of the $ SU_q(2) $ coherent state, namely, their contraction to 
%the coherent states for a quantum Heisenberg group, and incorporation of 
%the entanglement by employing a parity operator. To begin with we first 
study the contraction corresponding to the high spin limit. 
Such contraction of the algebras  $ \U $ and $ \A, $ that yields a quantum 
deformation of one dimensional Heisenberg algebra and group, respectively, 
has been discussed in \cite{CGST}. We apply the contraction of $ \U $ to our 
scheme. Consider the transformation of the generators of $\U:$
\beq
  A^{\dagger} = \frac{J_+}{ \sqrt{j} }, \qquad 
  A = \frac{J_-}{ \sqrt{j} }, \qquad H =  \frac{2}{j} J_0, 
      \qquad q = e^{w/j}.
  \label{contra-su2}
\eeq
The parameter $ w $ is assumed to be real number. 
Keeping $ w $ finite, we take the limit of $ j \rightarrow \infty. $ 
Then the commutation relations in (\ref{Uqdef}) yield
\beq
  [H, A] = [H, A^{\dagger} ] = 0, \qquad
  [ A, A^{\dagger} ] = \frac{\sinh wH}{w}.
  \label{Hq1def}
\eeq
The Hopf algebra mappings and $*$-involution are also contracted so that 
we have
\bea
  & & \Delta(H) = H \otimes 1 + 1 \otimes H, 
  \nn \\[3pt]
  & & \Delta(X) = X \otimes e^{wH/2} + e^{-wH/2} \otimes X, 
      \quad (X = A, \ A^{\dagger}),
     \label{Hq1Hopf} \\[3pt]
  & & \epsilon(X) = 0, \qquad S(X) = -X, \quad 
      (X = H,\ A,\ A^{\dagger}),
  \nn 
\eea
and
\beq
  (A^{\dagger})^* = A, \qquad A^* = A^{\dagger}, \qquad H^* = H.
  \label{starh1}
\eeq
In the classical limit of $ w \rightarrow 0 $ the relations (\ref{Hq1def}), 
(\ref{Hq1Hopf}) and (\ref{starh1}) reduce to their counterparts for the one 
dimensional Heisenberg algebra. Thus the Hopf $*$-algebra introduced above is 
a quantum deformation of the Heisenberg algebra, which we shall 
denote as $ U_q[h_1]. $

  We now turn to the contraction of the dual generators. 
The following transformations of generators of $ \A $ 
\beq
  \tilde{x} = \sqrt{j}\, x, \qquad
  \tilde{y} = \sqrt{j} \, y, \qquad
  \tilde{z} = \frac{j}{2} \, z 
  \label{contra-SU2}
\eeq
preserve the dual pairing between generators of 
$ \U $ and $ \A. $ Taking the limit $ j \rightarrow \infty, $ we obtain the 
commutation relations:
\beq
  [\tilde{x}, \tilde{y}] = 0, \qquad
  [\tilde{x}, \tilde{z}] = w \tilde{x}, \qquad 
  [\tilde{y}, \tilde{z}] = w \tilde{y}.
  \label{Hq1def-dual}
\eeq
At this stage we implement the contraction of the relation (\ref{uT2}) to 
extract the universal $\cal T$-matrix for the quantum Heisenberg algebra 
described above: 
\beq
  {\cal T} = \exp( \tilde{x} \otimes e^{-wH/2} A^{\dagger}) 
    \exp(\tilde{z} \otimes H) \exp( \tilde{y} \otimes e^{wH/2} A).
    \label{uTHq1}
\eeq
For the three dimensional representation $ \pi $ of the algebra $ U_q[h_1] $  
\beq
  \pi(A^{\dagger}) = \begin{pmatrix}
       0 & 0 & 0 \\
       0 & 0 & 1 \\
       0 & 0 & 0 \\
      \end{pmatrix},
     \quad
  \pi(A) = \begin{pmatrix}
       0 & 1 & 0 \\
       0 & 0 & 0 \\
       0 & 0 & 0 \\
      \end{pmatrix},
    \quad
  \pi(A^{\dagger}) = \begin{pmatrix}
       0 & 0 & 1 \\
       0 & 0 & 0 \\
       0 & 0 & 0 \\
      \end{pmatrix}
  \label{3Dreph1}
\eeq
the universal ${\cal T}$-matrix (\ref{uTHq1}) yields the matrix quantum group 
\beq
  (id \otimes \pi)({\cal T}) = 
  \begin{pmatrix}
       1 & y & z \\
       0 & 1 & x \\
       0 & 0 & 1 \\  
  \end{pmatrix}.
  \label{matrixH1}
\eeq
We observe that the matrix quantum group (\ref{matrixH1}) is identical to 
the quantum Heisenberg group $ H_q(1) $ given in \cite{CGST}. 
The Hopf algebra maps read \cite{CGST}:
\bea
  & & \Delta(\tilde{x}) = \tilde{x} \otimes 1 + 1 \otimes \tilde{x}, 
      \qquad
      \Delta(\tilde{y}) = \tilde{y} \otimes 1 +1 \otimes \tilde{y}, 
  \nn \\[3pt]
  & & \Delta(\tilde{z}) = \tilde{z} \otimes 1 + 1 \otimes \tilde{z} + \tilde{y} \otimes \tilde{x},
  \label{Hopf-Hq1} \\[3pt]
  & & \epsilon(a) = 0, \quad (a = \tilde{x}, \tilde{y}, \tilde{z}),
  \nn \\[3pt]
  & & S(\tilde{x}) = -\tilde{x}, \qquad S(\tilde{y}) = -\tilde{y}, \qquad 
      S(\tilde{z}) = -\tilde{z} + \tilde{x} \tilde{y}.
  \nn
\eea
The quantum group $H_q(1)$ is a Hopf $*$-algebra with the following 
involution map
\beq
  \tilde{x}^* = -\tilde{y}, \qquad \tilde{y}^* = -\tilde{x}, \qquad 
  \tilde{x}^* = -\tilde{z} + \tilde{x} \tilde{y}.
  \label{star-Hq1}
\eeq
It is obvious that $ {\cal T}^* {\cal T} = 1 \otimes 1. $

  Alternately, one can follow the prescription of Fronsdal and Galindo. 
Taking the basis $ E_{k \ell m} = (A^{\dagger})^k H^{\ell} A^m $ of 
$ U_q[h_1],$ 
we repeat the process in \S \ref{Pre}. The dual basis is determined to be
\beq
  e^{k\ell m} = 
  \frac{\tilde{x}^k}{k!} \frac{ (\tilde{z}- \hf (k-m) w)^{\ell} }{\ell !} 
  \frac{\tilde{y}^m}{m!}.
  \label{dualbasisHq1}
\eeq
The commutation relations (\ref{Hq1def-dual}), the Hopf structure 
(\ref{Hopf-Hq1}), and the universal ${\cal T}$-matrix (\ref{uTHq1}) are 
recovered confirming the validity of the contraction procedure. 

  To construct the coherent states for $ H_q(1) $ algebra, we first note 
its Fock space representation:
\bea
  & & A \ket{p\; n} = \sqrt{ n \frac{\sinh wp}{w} } \ket{p\; n-1}, 
  \nn \\[3pt]
  & & A^{\dagger} \ket{p; n} = \sqrt{ (n+1) \frac{\sinh wp}{w} } \ket{p\; n+1},
  \label{FockHq1} \\[3pt]
  & & H \ket{p; n} = p \ket{p; n}, \qquad p \in {\mathbb R}, \quad n = 0, 1, 2, \cdots.
  \nn
\eea
Parallel to the case of the ordinary bosons, the  coherent state for the  
$ H_q(1) $ algebra is constructed on the vacuum:
\beq
  \ket{\tilde{x}, \tilde{z}} = {\cal T} \ket{p\; 0} 
  = 
  e^{p\tilde{z}}\; \sum_{n=0}^{\infty} \left( \frac{ e^{2pw} - 1 }{2w} \right)^{n/2} 
  \frac{ \tilde{x}^n }{ \sqrt{n!} } \ket{p\;n}.
  \label{CS-Hq1}
\eeq
The factor $ e^{p \tilde{z}} $ normalizes the state appropriately.
This is verified by using the identity:
\[
  e^{p\tilde{z}^*} e^{p\tilde{z}} 
  = \exp\left( -\frac{ e^{2pw} - 1 }{2w} \tilde{x}^* \tilde{x}  \right).
\]
It is remarkable that  $ [\tilde{x}^*, \; \tilde{x}] = 0. $ Consequently, 
the K\"ahler geometry is almost trivial. 
However, the noncommutativity of $ \tilde{z} $ and $\tilde{x},\, \tilde{y} $ 
plays a crucial role regarding, for instance, the 
computations of expectation values and the resolution of unity. 
The resolution of unity here is much simpler than that of the $ SU_q(2)$ 
coherent states. Noting that 
\[
  e^{p\tilde{z}} e^{p\tilde{z}^*} = 
  \exp\left( -\frac{ 1 - e^{-2pw}  }{2w} \tilde{x}^* \tilde{x}  \right),
\]
the invariant integration is reduced to the usual integration on 
the complex plane by regarding $\tilde{x}$ as a complex variable
and $ \tilde{x}^* $ as its conjugate. Setting $ \tilde{x} = r e^{i\theta}, $  
it is easy to verify that
\beq
 \int d\mu \ket{\tilde{x}, \tilde{z}} \bra{\tilde{x}, \tilde{z}} = 1, \qquad
 d\mu = \frac{ 1 - e^{-2pw}  }{2\pi w} r dr d\theta.
 \label{DU-Hq1}
\eeq

%%%%%%%%%%%%%%%%%%%%%%%%%%%%%%%%%%%%%%%%%%%%%%%%%%%%%%%%%%%%
%

%%%%%%%%%%%%%%%%%%%%%%%%%%%%%%%%%%%%%%%%%%%%%%%%%%%%%%%%%
%
%  Concluding remarks
%
%%%%%%%%%%%%%%%%%%%%%%%%%%%%%%%%%%%%%%%%%%%%%%%%%%%%%%%%%
%
\section{Concluding remarks}
\setcounter{equation}{0}
\label{CRem}

 We have investigated the $ SU_q(2) $ coherent states 
in detail. It was shown that properties analogous to the classical 
$ SU(2) $ coherent states also hold for 
its quantum group counterpart. A characteristic feature of $ SU_q(2) $ 
coherent states is known to be the noncommutativity of 
the variable parametrizing the states. Thanks to this 
fact, we obtained a natural description 
of the $q$-sphere in complex coordinates. Our description of the differential 
calculus on the complexified $q$-sphere may provide essential tools for 
constructing the path integrals on it. In addition, this may pave a way for 
studying the K\"ahler structure on the $q$-sphere. Probably one can generalize 
this to other noncommutative analogues of K\"ahlerian manifolds. 
Furthermore, similarity of representation theory between 
$ su(2) $ and the Lie superalgebra $ osp(1/2) $ encourages us 
to study coherent state for the quantum supergroup $ OSp_q(1/2) $\cite{KR}. 
The universal $\cal T$-matrix for $ OSp_q(1/2) $ has been obtained \cite{ACNS} 
and the finite dimensional representations of $ OSp_q(1/2) $ are well studied 
\cite{KR,Zou,ACNS,ACNS2}. 
We are ready to study coherent states for $ OSp_q(1/2). $ 
Once the  $ OSp_q(1/2) $ coherent states are obtained, 
they will give a complex description 
of $q$-supersphere introduced in \cite{AC}. 
Along the line similar to the classical case \cite{Gra}, 
one will be able to discuss noncommutative version of super-K\"ahler geometry 
using the $ OSp_q(1/2) $ coherent states. 
This work is in progress. 

  We have made important observations such as the resolution of unity
in the context of the $ SU_q(2) $ coherent states. 
As is well known, many applications of 
coherent states stem from this property. 
We thus believe the coherent states discussed in this 
paper has potential of various applications in 
physics or mathematics where noncommutativity plays 
certain roles. One of such possibilities may be 
coherent state description of operators. As shown in 
\cite{Lieb}, operators appearing in the analysis of 
spin system are described by spherical harmonics. 
It may be expected that consideration of operators 
consisting of the generators of $ \U $ leads us 
to define noncommutative version of spherical harmonics. In another development
\cite{BCHR} it has been observed that the effective $su_{q}(2)$ hamiltonians
successfully reproduce the ground state properties and the spectrum of 
different interacting fermion-boson dynamical nuclear systems. The bosonic 
part of the interactions can be effectively embedded as an appropriate 
$q$-deformation of the fermionic $su(2)$ algebra. The resolution of unity via 
the coherent states obtained here may be useful in obtaining suitable matrix 
elements in these models.

%%%%%%%%%%%%%%%%%%%%%%%%%%%%%%%%%%%%%%%%%%%%%%%%%%%%%%%%%
%
%  Appendix 
%
%%%%%%%%%%%%%%%%%%%%%%%%%%%%%%%%%%%%%%%%%%%%%%%%%%%%%%%%%
%
\section*{Appendix: $q$-Analysis}
\setcounter{equation}{0}
\renewcommand{\theequation}{A.\arabic{equation}}

  Formulae of $q$-analysis used in this paper are summarized in 
this appendix.  We follow the notations and conventions of \cite{GR}. 
\begin{enumerate}
  \item $q$-shifted factorial
  \bea
    & & (a;q)_n = \left\{
          \begin{array}{lcl}
            1, & & n = 0 \\[5pt]
            {\displaystyle \prod_{k=0}^{n-1} (1 -a q^k)}, & & n \in {\mathbb N}
          \end{array}
        \right.
        \label{qshifted1} \\[5pt]
    & & (a;q)_{\infty} = {\displaystyle \prod_{k=0}^{\infty} (1 -a q^k)}.
        \label{qshifted2}
  \eea
  \item useful identities
  \bea
    & & (a;q)_n = \frac{ (a;q)_{\infty} }{ (aq^n;q)_{\infty} }, 
        \qquad
        \frac{1}{(a;q)_n} = \frac{ (a q^n;q)_{\infty} }{ (a;q)_{\infty} },
        \label{qid1} \\[5pt]
    & & \frac{ (q;q)_m }{ (q;q)_{m-n} } = (-1)^n q^{mn-n(n-1)/2} (q^{-m};q)_n, 
        \label{qid2} \\[5pt]
    & & \bn1{m}{k} = (-1)^k q^{(m+1)k} \frac{ (q^{-2m};q^2)_k }{ (q^2;q^2)_k }.
        \label{bn-qfac}
  \eea
  \item basic hypergeometric series
  \bea
    & & \BH{r}{s}\BHA{a_1,a_2,\cdots,a_r}{b_1,\cdots, b_s}{q}{z} 
      \nn \\[5pt]
    & & \qquad = 
        \sum_{n=0}^{\infty} \frac{(a_1;q)_n (a_2;q)_n \cdots (a_r;q)_n}
                             {(q;q)_n (b_1;q)_n \cdots (b_s;q)_n} 
        \{ (-1)^n q^{n(n-1)/2} \}^{1+s-r} z^n.
        \label{BHSdef}
  \eea
  \item $q$-binomial theorem
  \beq
    \BH{1}{0}\BHA{a}{-}{q}{z} 
    = \sum_{n=0}^{\infty} \frac{ (a;q)_n }{ (q;q)_n } z^n
    = \frac{ (az;q)_{\infty} }{ (z;q)_{\infty} }.
    \label{q-bThm}
  \eeq
  \item special case of $ \BH{1}{1} $
  \beq
    \BH{1}{1}\BHA{a}{c}{q}{\frac{c}{a}} = \frac{ (c/a;q)_{\infty} }{ (c;q)_{\infty} }.
    \label{11trans}
  \eeq
  \item $q$-derivative and $q$-integral
  \bea
    & & D_q f(x) = \frac{ f(x) - f(xq) }{ (1-q)x }, \label{q-der} \\[4pt]
    & & \int_0^x f(t) d_qt = (1-q)x \sum_{k=0}^{\infty} f(xq^k) q^k,
        \label{q-int} \\[4pt]
    & & \int_0^x D_q f(t) d_qt = D_q \int_0^x f(t) d_qt = f(x).
        \label{der-int}
  \eea
  Leibniz rule
  \beq
    D_q f(x) g(x) = (D_q f(x)) g(xq) + f(x) D_q g(x). \label{q-Leib}
  \eeq
  Integral by parts
  \bea
   & & 
    \int_0^x (D_q f(t)) g(tq) d_qt = f(x) g(x) - \int_0^x f(t) D_q g(t) d_qt,
    \label{intpart} \\[5pt]
   & & 
    \int_0^x (D_q f(t)) g(t) d_qt = f(x) g(x) - \int_0^x f(tq) D_q g(t) d_qt.
    \label{intpart2}
  \eea
  Some formulae
  \bea
    & & D_q \, x^n = (n)_q \, x^{n-1}, \qquad
        D_q (xq^{-1};q)_n = -q^{-1} (n)_q \, (x;q)_{n-1},
        \label{der-form} \\[5pt]
    & & D_q (q^{-a+k} x; q)_a = -q^{-a+k} (a)_q (q^{-a+k+1}x;q)_{a-1},
        \label{der-form2} \\[5pt]
    & & \int_0^x t^n d_qt = \frac{ x^{n+1} }{ (n+1)_q }.
        \label{int-form}
  \eea
\end{enumerate}

\vspace{8mm}\noindent
\textit{Note added} $-$ After submission of our paper, we were informed that 
\v{S}koda introduced coherent states for Hopf algebras 
based on quantum line bundles \cite{Sko}. As an example, the coherent state 
of the $SU_q(2)$ algebra enjoying a resolution of unity was discussed. 
We thank Zoran \v{S}koda for sending us his work.

%%%%%%%%%%%%%%%%%%%%%%%%%%%%%%%%%%%%%%%%%%%%%%%%%%%%%%%%
%
%  REFERENCES
%
%%%%%%%%%%%%%%%%%%%%%%%%%%%%%%%%%%%%%%%%%%%%%%%%%%%%%%%%%
%

\end{document}